\documentclass[11pt]{amsart}
\usepackage{amsfonts,amssymb,amscd,amsmath,enumerate,verbatim,calc}

\def\NZQ{\mathbb}               
\def\NN{{\NZQ N}}
\def\QQ{{\NZQ Q}}
\def\ZZ{{\NZQ Z}}

\def\frk{\mathfrak}               

\def\pp{{\frk p}}

\def\qq{{\frk q}}

\def\mm{{\frk m}}

\def\Phi{{\frk N}}
\def\Aa{{\frk A}}
\def\opn#1#2{\def#1{\operatorname{#2}}} 
\opn\chara{char} \opn\length{\ell} \opn\pd{pd} \opn\rk{rk}
\opn\projdim{proj\,dim} \opn\injdim{inj\,dim} \opn\rank{rank}
\opn\depth{depth} \opn\grade{grade} \opn\height{height}
\opn\embdim{emb\,dim} \opn\codim{codim}

\opn\Tr{Tr} \opn\bigrank{big\,rank}
\opn\superheight{superheight}\opn\lcm{lcm}
\opn\trdeg{tr\,deg}
\opn\reg{reg} \opn\lreg{lreg} \opn\ini{in} \opn\lpd{lpd}
\opn\size{size}\opn{\mult}{mult}
\opn\div{div} \opn\Div{Div} \opn\cl{cl} \opn\Cl{Cl}
\opn\Spec{Spec} \opn\Supp{Supp} \opn\supp{supp} \opn\Sing{Sing}
\opn\Ass{Ass} \opn\Min{Min}
\opn\Ann{Ann} \opn\Rad{Rad} \opn\Soc{Soc}
\opn\Syz{Syz} \opn\Im{Im} \opn\Ker{Ker} \opn\Coker{Coker}
\opn\Am{Am} \opn\Hom{Hom} \opn\Tor{Tor} \opn\Ext{Ext}
\opn\End{End} \opn\Aut{Aut} \opn\id{id}

\opn\nat{nat}
\opn\pff{pf}
\opn\Pf{Pf} \opn\GL{GL} \opn\SL{SL} \opn\mod{mod} \opn\ord{ord}
\opn\Gin{Gin}
\opn\Hilb{Hilb}\opn\adeg{adeg}\opn\std{std}\opn\ip{infpt}
\opn\Pol{Pol}
\opn\sat{sat}
\opn\Var{Var}
\opn \ann{ann}
\opn\aff{aff} \opn\con{conv} \opn\relint{relint} \opn\st{st}
\opn\lk{lk} \opn\cn{cn} \opn\core{core} \opn\vol{vol}
\opn\link{link} \opn\star{star}
\opn\gr{gr}
\def\Rees{{\mathcal R}}

\def\Fc{{\mathcal F}}
\def\pot#1#2{#1[\kern-0.28ex[#2]\kern-0.28ex]}

\opn\dirlim{\underrightarrow{\lim}}
\opn\inivlim{\underleftarrow{\lim}}
\let\tensor=\otimes

\let\Union=\bigcup
\let\Sect=\bigcap
\let\Dirsum=\bigoplus

\let\to=\rightarrow
\let\To=\longrightarrow
\def\Implies{\ifmmode\Longrightarrow \else
        \unskip${}\Longrightarrow{}$\ignorespaces\fi}
\def\implies{\ifmmode\Rightarrow \else
        \unskip${}\Rightarrow{}$\ignorespaces\fi}
\def\iff{\ifmmode\Longleftrightarrow \else
        \unskip${}\Longleftrightarrow{}$\ignorespaces\fi}

\let\:=\colon
\theoremstyle{plain}
\newtheorem{Theorem}{Theorem}[section]
\newtheorem{Lemma}[Theorem]{Lemma}
\newtheorem{Corollary}[Theorem]{Corollary}
\newtheorem{Proposition}[Theorem]{Proposition}
\newtheorem{Remark}[Theorem]{Remark}

\newtheorem{Example}[Theorem]{Example}

\newtheorem{Definition}[Theorem]{Definition}

\let\epsilon\varepsilon
\let\phi=\varphi
\let\kappa=\varkappa
\let\ol=\overline
\def\qed{\ifhmode\textqed\fi
      \ifmmode\ifinner\quad\qedsymbol\else\dispqed\fi\fi}
\def\textqed{\unskip\nobreak\penalty50
       \hskip2em\hbox{}\nobreak\hfil\qedsymbol
       \parfillskip=0pt \finalhyphendemerits=0}
\def\dispqed{\rlap{\qquad\qedsymbol}}

\opn\dis{dis}
\def\pnt{{\raise0.5mm\hbox{\large\bf.}}}

\opn\Lex{Lex}

\begin{document}

\title{Hilbert polynomials and powers of ideals}

\author{J\"urgen Herzog, Tony J. Puthenpurakal and Jugal K. Verma}

\address{J\"urgen Herzog, Fachbereich Mathematik und
Informatik, Universit\"at Duisburg-Essen, Campus Essen, 45117
Essen, Germany} \email{juergen.herzog@uni-essen.de}

\address{Tony J. Puthenpurakal, Department of Mathematics, Indian Institute
of Technology, Mumbai,
India}\email{tputhen@math.iitb.ac.in}

\address{Jugal K. Verma, Department of Mathematics, Indian Institute
of Technology, Mumbai,
India}\email{jkv@math.iitb.ac.in}

\thanks{Jugal Verma and Tony J. Puthenpurakal  thank Universt{\"a}t
Duisburg-Essen
  for hospitality during Sept-Dec 2006.
The authors thanks DFG for financial support, which made this visit possible}

\begin{abstract}
The growth of Hilbert coefficients for powers of ideals are studied. For a graded ideal $I$ in the polynomial ring $S=K[x_1,\ldots,x_n]$ and a finitely generated graded $S$-module, the Hilbert coefficients $e_i(M/I^kM)$ are polynomial functions.    Given two families of graded ideals $(I_k)_{k\geq 0}$ and $(J_k)_{k\geq 0}$ with $J_k\subset I_k$ for all $k$ with the property that $J_kJ_\ell\subset J_{k+\ell}$ and $I_kI_\ell\subset I_{k+\ell}$ for all $k$ and $\ell$, and such that the algebras $A=\Dirsum_{k\geq 0}J_k$ and $B=\Dirsum_{k\geq 0}I_k$ are finitely generated, we show  the function $k \mapsto_0(I_k/J_k)$ is of quasi-polynomial type, say given by the polynomials $P_0,\ldots, P_{g-1}$. If
$J_k = J^k$ for all $k$ then we show that all the $P_i$ have the same degree and the same leading coefficient.
As one of the  applications it is shown that $\lim_{k\to \infty}\length(\Gamma_\mm(S/I^k))/k^n \in \mathbb{Q}.$
We also study analogous statements in the local case.
\end{abstract}
\maketitle
\section*{Introduction}

This paper is inspired by a result of Cutkosky, Ha, Srinivasan and Theodorescu. In their paper \cite{CHT} they showed that for a graded ideal $I$ in a polynomial ring $K[x_1,\ldots,x_n]$ over a field $K$ of characteristic  $0$, the limit $\lim_{k\to \infty} \length((I^k)^{\sat}/I^k)/k^n$ exists, but need not to be a rational number. Here $I^{\sat}$ denotes the saturation of $I$ and $\length(M)$ the length of a module.   Our primary question then was whether this limit is a rational number when $I$ is a monomial ideal. This is indeed the case,  as we show in Corollary \ref{conclusion}.

More generally,  one may consider the $j$th symbolic power of an ideal $I$ with respect to an ideal $J$, namely the ideal $I_k(J)=I^k\: J^\infty$. The case of saturated powers is the special case in which $J$ is the graded maximal ideal of $S$.  Again, if $I$ and $J$ are monomials we can show that $\lim_{k\to \infty} e_0(I_k(J)/I^k)/k^{n-d}$ exists and is a rational number. Of course, $I_k(J)/I^k$ is not a finite length module in general, so that in this limit formula we had to replace the length of  $I_k(J)/I^k$ by the multiplicity  $e_0(I_k(J)/I^k)$ of this module. Moreover in the denominator we can replace  $k^n$  by $k^{n-d}$, where $d$ is the limit dimension, that is, the dimension of  $I_k(J)/I^k$ for $k\gg 0$, as it can be shown that  $\dim I_k(J)/I^k$ is constant for large $k$.
For the proof of this result we use the fact  that the graded $S$-algebra $\Dirsum_{k\geq 0}I_k(J)$ is finitely generated if $I$ and $J$ are monomial ideals, as shown in \cite[3.2]{HHT}.

The result described in the preceding paragraph is again a special case of the following more general situation: given  two families $(I_k)_{k\geq 0}$ and $(J_k)_{k\geq 0}$ of graded ideals of $S$ with the property that $J_kJ_\ell\subset J_{k+\ell}$ and $I_kI_\ell\subset I_{k+\ell}$ for all $k$ and $\ell$, and such that the algebras $A=\Dirsum_{k\geq 0}J_k$ and $B=\Dirsum_{k\geq 0}I_k$ are finitely generated. We also assume that $J_k\subset I_k$ for all $k$ and that $I_0=J_0=S$.

What can be said about  $e_0(I_k/J_k)$ as a function of $k$? We first notice in Proposition \ref{quasi} is that $e_0(I_k/J_k)$ is a quasi-polynomial for $k\gg 0$. Indeed, this can be deduced from the fact that $e_0(I^kM/J^kM)$ is quasi-polynomial for $k\gg 0$ for any finitely generated graded $S$-module $M$ and graded ideal $J\subset I\subset S$, see Proposition \ref{limit}. For the proof of this proposition we use Theorem \ref{graded}, where we consider a finitely generated graded $S$-module $N$. The Hilbert polynomial  $P_N(x)$ of $N$  can be
written in the form
\begin{eqnarray}
P_N(x)=\sum_{i=0}^{d}(-1)^ie_{i}(N)\binom{x+d-i}{d-i}
\end{eqnarray}
with integer coefficients $e_i(N)$, called the {\em Hilbert
coefficients} of $N$. Now Theorem \ref{graded} states the following: Let $M$ be a graded $S$-module, and $I\subset S$ a graded ideal with $\dim M/IM=d$. Then for all $i=0,\ldots, d$, the
Hilbert coefficient $e_{i}(M/I^kM)$  as a function of $k$ is of
polynomial type of degree $\leq n-d+i$. The assertion that $e_{i}(M/I^kM)$ is a polynomial follows easily from  \cite[4.3]{HoangT} (also see proof of \cite[Theorem 1]{Ko}). However our
degree bound is new.

In Section 5 we give a different proof of Proposition \ref{limit}, not referring to Theorem \ref{graded}. This alternative proof is also valid for ideals $J\subset I$ in a local ring. In fact, we do not have a local version of Theorem \ref{graded}. In other words, we do not know whether  for an ideal $I$ in a local ring $(R,\mm)$ the Hilbert coefficients $e_i(R/I^k)$ are polynomial functions for large $k$. This is well-known and easy to prove for  $i=0$, but seems to be unknown for $i>1$. Actually we expect that $e_i(R/I^k)$ may not be of polynomial type for $i>1$.

For any quasi-polynomial numerical function like $e_0(I_k/J_k)$  (say of period $g$ and given by polynomials $P_0, \ldots, P_{g-1}$) a natural question is:
\begin{equation*}
 \text{Do the polynomials $P_0, \ldots, P_{g-1}$ have the same degree and the same leading term?}  \tag{$*$}
\end{equation*}
To put this question in the right historical framework let us consider the following:
Let $P$ be a $d$-dimensional rational covex polytope in $\mathbb{R}^n$ and let $E(P,k) = \sharp \ kP \cap \mathbb{Z}^n $ be the Ehrhart
function. It is well-known that   $E(P,k)$  is of quasi-polynomial type cf. \cite[6.3.11]{BH}, say given by
$Q_0, \ldots, Q_{g-1}$. We state a
simplified version of a conjecture due to Ehrhart: \textit{If the affine span of every $d-1$ dimensional face of $P$ contains a point with integer co-ordinates then
$Q_0, \ldots, Q_{g-1}$ have the same degree and the same leading term.}

A more general version of this was proved independently by  McMullen (see \cite{M}) and Stanley
(\cite{S}, Theorem 2.8). For a recent proof see \cite[5]{BI}.

In our case, since we know that $e_0(I_k/J_k)$ is of quasi-polynomial type, we also would like to know whether the polynomials describing this quasi-polynomial function  are all of same degree and all have the same leading term. Then this would imply at once  that $\lim_{k \rightarrow \infty} e_0(I_k/J_k)/k^{n-d}$ exists and is a rational number.

The right set-up to study  these questions is to consider the case of homogeneneous inclusion of Noetherian algebras
$A = \bigoplus_{k \geq 0}A_k \subseteq B = \bigoplus_{k \geq 0}B_k$, with $A_0 = B_0$. We assume $A$ is standard
graded, while $B$ need not be. Furthemore $A_0$ is either the polynomial ring  $S = K[x_1,\ldots,x_n]$ or  a local Noetherian ring $(R,\mm)$
with residue field $K = R/\mm$. In the case when $A_0 = S$ we also assume that $A$, $B$ are bigraded.
In this general set-up one can similarly prove that $k \mapsto e_0(B_k/A_k)$ is of quasi-polynomial-type (in the local case we take
Hilbert-Samuel multiplicity with respect to $\mm$). We consider the question ($*$) for this function.

Let $M = \bigoplus_{k \geq 0}M_k$ be a graded but not necessarily finitely generated $A$-module (it is bigraded if $A$ is). We say $M$ is a quasi-finite $A$-module if  $H^0_{A_+}(M)_i=0$ for $i\gg 0$.  Here  $H^0_{A_+}(M)$ denotes the $0$th local cohomology of $M$ with respect to the ideal $A_+=\Dirsum_{k>0}A_k$.
Our main result (see Theorem \ref{tony} and Theorem \ref{MTinclusion})  is that when $B/A$ is quasi-finite $A$-module, then indeed the function $k \mapsto e(B_k/A_k)$ satisfies ($*$). The module $B/A$ is quasi-finite under the mild hypothesis that
$\grade(A_1B, B) >0$. This is trivially satisfied  in the cases when $A$ is the Rees algebra $\Rees(I)=\bigoplus_{n \geq 0} I^nt^n$ of an ideal $I$ and $B = \bigoplus_{k \geq 0}I_k$, where  $(I_k)_k\geq 0$  is a family of  ideals such that $I^k \subseteq I_k$ for all $k$.  In the case when $A_0 = R$ is local,
$A = \Rees(J)$ and $B =\Rees(\Fc)$, then note that $\grade(A_1B, B) >0$ if $\grade(J,R) > 0$.

The main technical tool required is that if $M$ is a quasi-finite $A$-module and if $K$ is uncountable, then $M$ has a filter-regular element $x \in A_1$. Filter-regular elements in the case when $M$ is finitely generated have a lot of  applications in the study of blow-up algebras \cite{t2}. An early version of quasi-finite modules is in \cite[4.7]{TJP}. However it was somewhat  unexpected to us that the existence of filter-regular elements for quasi-finite modules have implications in question ($*$) above.

In section four, we specialise the above set-up to the  case  where $A=\Rees(J) \subset B=\Rees(I)$ for ideals $J \subseteq I$ of a $d$-dimensional local ring $(R,\mm)$
and $\ell(I/J) < \infty.$ Amao \cite{a} showed that the numerical function $n \mapsto H(I/J,n):=\ell(I^n/J^n)$  is a polynomial function of degree atmost $d.$ We call $H(I/J,n)$ the Rees function of the pair $(J,I)$ and the corresponding polynomial $P(I/J,n),$ the Rees polynomial of the
pair $(J,I).$   Our objective is to identify the degree of the Rees polynomial.  Rees \cite{r} showed that if $(R, \mm)$ is quasi-unmixed,
then $J$ is a reduction of $I$ if and only if $\deg P(I/J,n) < d.$  When $I$ and $J$ are $\mm$-primary,
identification of the degree of $P(I/J,n)$ amounts to  knowing the number of Hilbert-Samuel coefficients  of $I$ and $J$ which coincide.  This problem was considered in  \cite{s} by K. Shah. He introduced
coefficient  ideals of an $\mm$-primary ideal which help in identification of the degree of the
Rees polynomial. We show in Theorem  \ref{coefficient} that coefficient ideals exist even in the general case. Computation of coefficient ideals seems to be a difficult problem at present. We show that
when $J$ is a reduction of $I$ then the degree of the Rees polynomial of the pair $(J,I)$ is
$\dim \Rees(I)/\Rees(J)-1.$  As a consequence it follows that if $R$ is quasi-unmixed, $J \neq I$ and $\Rees(J)$  satisfies Serre's condition $S_2,$ then $\deg P(I/J,n)=d-1.$

Recall that an ideal is called normal if all its powers are integrally closed.
It is rather surprising that for distinct normal ideals $J \subset I,$ the problem of
determination of the degree  of the Rees polynomial of the pair $(J,I)$ can be completely solved.
Let $J \subset I$ be ideals of an analytically unramified $d$-dimensional local ring $R$ such that $\ell(I/J)$ is finite and $J$ is not a reduction of $I.$
 In Theorem  \ref{normalRees}
we show that $f(n):=\ell(\ol{I^n}/\ol{J^n})$ is  a polynomial function of degree $d.$
It follows that if $J$ and $I$ are normal  then $\deg P(I/J,n)=d.$

In  Section 5,  $(R, \mm)$ is a Noetherian local ring, $I$ an ideal in $R$, $M$  a finitely generated $R$-module and $J$ a reduction of $I$ with respect to $M$, i.e., $I^{m+1}M = JI^mM$ for some $m$.
If $\dim M = r$, then for $i = 0,1,\ldots, r$ let  $e_i(M)$ be the $i$th Hilbert coefficient of $M$ (with respect to $\mm$). In other words,
\[
\ell(M/\mm^{n+1}M) = \sum_{i = 0}^{r}(-1)^ie_i(M)\binom{n+ r -i}{r-i} \quad \text{for all} \ n \gg 0.
\]

By a result due to Brodmann \cite{Br}, $\Ass M/I^kM$ is stable for $k \gg 0$. It follows that
$\dim M/I^kM$ is stable for $k \gg 0$. We denote this stable value by $\theta_I(M)$.
It is clear that $ \bigoplus_{k \geq 0} I^kM/I^{k+1}M $ and $\bigoplus_{k \geq 0} I^kM  $ are finitely generated
modules over $\Rees(I).$ Since
$ \bigoplus_{k \geq 0} I^kM/J^{k}M$ is a
 finitely generated module over $ \Rees(J),$   $\dim I^k M$, $\dim I^kM/I^{k+1}M$
and $\dim I^kM/J^kM$ stabilize  for $k \gg 0$. We denote these stable values by  $\alpha_I(M), \beta_I(M)$ and $\gamma^{I}_{J}(M)$, respectively.

In view of Theorem \ref{graded} one may ask the
following question: Is it true that  the numerical functions
\begin{enumerate}[\rm (1)]
\item $e_i(M/I^{k}M)$,  $i = 0,\ldots, \theta_I(M)$,
\item $e_i(I^kM/I^{k+1}M)$, $i = 0,\ldots, \alpha_I(M)$,
\item $e_i(I^kM)$,   $i = 0,\ldots, \beta_I(M)$, and
\item $e_i(I^k M/J^{k}M)$, $i = 0,\ldots, \gamma^{I}_{J}(M)$
\end{enumerate}
are polynomial functions  in $k$ ?

The graded versions of  Question (1),(2) and (3) are equivalent, while  (4) follows from (2). However in the local case all these questions are different. In Section 5, we settle Question (2), (3) and (4) (see Corollary \ref{locT}). As we mentioned already, the answer to Question (1) is still open. It can be easily checked that if $x^*$, the initial form of $x$ in the associated  graded $\gr_{\mm}(R)$ of $R$, is  $\gr_{\mm}(R)$-regular then $e_i(R/(x)^k)$ is a polynomial in $k$. We do not know of any other cases where this question  has a positive answer.

\section{The Hilbert coefficients of $M/I^kM$ as a function of $k$}
Let $K$ be a field, $S=K[x_1,\ldots, x_n]$ the polynomial ring in
$n$ variables,  and let $N$ be any graded $S$-module of dimension
$d$. Then for $i\gg 0$, the numerical function $H(N,i)=\sum_{j\leq
i}\dim_KN_j$ is a polynomial function of degree $d$, see \cite[4.1.6]{BH}.
In other words, there exists a polynomial $P_N(x)\in \ZZ[x]$ such
that
\[
H(N,i)=P_N(i) \quad \text{for all}\quad i\gg 0.
\]
The polynomial $P_N(x)$ is called the {\em Hilbert polynomial} of $N$.
It can be
written in the form

\begin{eqnarray}
\label{formula}
P_N(x)=\sum_{i=0}^{d}(-1)^ie_{i}(N)\binom{x+d-i}{d-i}
\end{eqnarray}
with integer coefficients $e_i(N)$, called the {\em Hilbert
coefficients} of $N$.

Let $I\subset S$ a graded ideal and $M$ a graded $S$-module with
$\dim M/IM=d$. We are interested in the Hilbert coefficients
$e_i(M/I^kM)$ of $M/I^kM$ as functions in $k$.

We say that a function $f\:\NN \to \QQ$ is of polynomial type of
degree $d$, if there exists a polynomial $h\in \QQ[x]$ of degree
$d$ such that $f(k)=h(k)$ for all $k\gg 0$.

As the main result of this section we prove

\begin{Theorem}
\label{graded}
 Let $M$ be a graded $S$-module, and $I\subset S$ a
graded ideal with $\dim M/IM=d$. Then for all $i=0,\ldots, d$, the
Hilbert coefficient $e_{i}(M/I^kM)$  as a function of $k$ is of
polynomial type of degree $\leq n-d+i$.
\end{Theorem}

\begin{proof}
For the proof  we proceed by induction on $d$. If $d=0$, then
$P_{M/I^kM}(x)=\length(M/I^kM)$, and hence the only Hilbert
coefficient of $M/I^kM$ is $e_0(M/I^kM)=\length(M/I^kM)$. This
length, as a function of $k$, is known to be of polynomial type of
degree $\leq n$, see \cite[4.1.6]{BH}.

Now let $d>0$. By a theorem of Brodmann \cite{Br} the set
$\Union_k\Ass(M/I^kM)$ is finite. Hence if $|K|=\infty$, as we may
assume, there exists a linear form $y\in S$ such that the kernel
of the multiplication map $y\: (M/I^kM)(-1)\to M/I^kM$ is of
finite length for all $k$. Set $N=M/y M$, then $N/I^kN=
(M/I^kM)/y(M/I^kM)(-1)$, and it follows that
$$P_{N/I^kN}(j)=P_{M/I^kM}(j)-P_{M/I^kM}(j-1) \quad\text{for all}\quad j.$$
From this we conclude that $e_i(M/I^kM)=e_i(N/I^kN)$ for
$i=0,\ldots, d-1$. By induction hypothesis $e_i(N/I^kN)$ is of
polynomial type  of degree $\leq
 (n-1)-(d-1)+i=n-d+i$.  Thus it remains to prove that
 $e_{d}(M/I^kM)$ is  of  polynomial type  of degree $\leq
 n$.

We note that for all integers $a$ one has
\[
e_{d}(M/I^kM)=P_{M/I^kM}(ak)-\sum_{i=0}^{d-1}(-1)^ie_{i}(M/I^kM)\binom{ak+d-i}{d-i},
\]
and that each of the summands $e_{i}(M/I^kM)\binom{ak+d-i}{d-i}$ is
of  polynomial type   of degree $\leq n$. Thus it suffices to show
that $P_{M/I^kM}(ak)$ is of polynomial type  of degree $\leq n$.

To this end, we consider the sequence
\[
0\To I^kM\To M\To M/I^kM\To 0,
\]
from which we deduce that
$$P_{M/I^kM}(ak)=P_{M}(ak)-P_{I^kM}(ak).$$ Since $P_{M}(x)$ is a
polynomial of degree $\leq n$ it follows that $P_{M}(ak)$ is a
polynomial in $k$ of degree $\leq n$. Thus the assertion follows
once it is shown that for suitable integer $a>0$, $P_{I^kM}(ak)$
is of polynomial type of degree $\leq n$.

We consider the Hilbert series
\[
\Hilb_{I^kM}(t)=\sum_{i\geq 0}
H(I^kM,i)t^i=\frac{Q_k(t)}{(1-t)^{n+1}}
\]
of $I^kM$. Here $Q_k(t)$ is a polynomial. Indeed let
$\beta_{ij}=\dim_K\Tor_i^S(I^kM,K)_j$ be the graded Betti-numbers
of $I^kM$. Then
$$
Q_k(t)=\sum_{i=0}^n(-1)^i\beta_{ij}t^j,
$$
see \cite[4.1.13]{BH}. It follows that $\deg Q_k(t)\leq \reg(I^kM)+n$.
As shown in \cite[2.4]{CHT} and \cite[1]{Ko} one has $\reg(I^kM)\leq ek+f$
for all $k$ and certain integers $e$ and $f$ with $e>0$.
Therefore, $\deg \Hilb_{I^kM}(t)= ek+f-1<ak$ for all $k$ and any
integer $a>\min\{e,e+f-1\}$.
By \cite[4.1.12]{BH} we have $H(I^kM,i)=P_{I^kM}(i)$ for all $i>\deg
\Hilb(t)$. In particular,
$H(I^kM,ak)=P_{I^kM}(ak) \quad \text{for all} \quad k.$
Thus it remains  to show that  $H(I^kM, ak)$ is of polynomial type
of degree $\leq n$ for a suitable $a$.

Let $W=M\tensor_SS[x_{n+1}]$; then $(I^kW)_j=\Dirsum_{i\leq
j}(I^kM)_ix_{n+1}^{j-i}$, so that $H(I^kM, ak)=\dim_K(I^kW)_{ak}$
for all $k$. Hence we must show that $\dim_K(I^kW)_{ak}$ is of
polynomial type of degree $\leq n$ for a suitable $a$. Since
$I^kS[x_{n+1}]W=I^kW$, we may replace $I^k$ by its extension
$I^kS[x_{n+1}]=(IS[x_{n+1}])^k$ which, by an abuse of notation,
we again denote by $I^k$

Let $c_k$ be the highest degree of a generator of $I^kW$, $c$ the
highest degree of a generator of $I$ and $d$ the highest degree of
a generator of $W$. Then one has $c_k\leq kc+d$ for all $k$. Let
$J=I_{\geq c}=\Dirsum_{j\geq c} I_j$. We choose $a>c+d$; then
\begin{eqnarray*}
\dim_K(J^kW)_{ak}=\dim_K((J^kW)_{\geq
c_k})_{ak}=\dim_K((I^kW)_{\geq c_k})_{ak} &=&\dim_K(I^kW)_{ak}.
\end{eqnarray*}
Thus we may replace $I$ by $J$. Since $J$ is generated in one
degree, we may as well assume that $I$ is generated in one degree,
say $c$. In this case the Rees algebra ${\mathcal
R}(I)=\Dirsum_{j\geq 0}I^j$ is naturally standard bigraded, the
bigraded components being ${\mathcal R}(I)_{(g,h)}=(I^h)_{g+ch}$.
Moreover, $V=\Dirsum_{k\geq 0}I^kW$ is a finitely generated
bigraded ${\mathcal R}(I)$-module with $V_{(g,h)}=(I^hW)_{g+ch}$.

If $R$ is any standard  bigraded $K$-algebra, and $e$ and $f$ are
two natural numbers. Then the subalgebra $R_\Delta=\Dirsum_{k\geq
0}R_{(ke,kf)}$ is  called the $(e,f)$-diagonal of $R$. It is
easily seen that any such diagonal is a standard graded
$K$-algebra. Moreover, if $N$ is a finitely generated bigraded $R$-module, then the
diagonal submodule $N_\Delta=\Dirsum_{k\geq 0}N_{({ek,fk})}$ is a
finitely generated graded $R_\Delta$-module.

We now consider the $(a-c,1)$-diagonal $\Delta$ of ${\mathcal
R}(I)$. Then $V_\Delta=\Dirsum_{k\geq 0}(I^kW)_{ak}$ is a finitely
generated graded $A$-module where $A$ is the   standard graded
$K$-algebra with $A_k=(I^k)_{ak}$ for all $k$. It follows that
$\dim_K(I^kW)_{ak}=\dim_K(V_{\Delta})_k$ is of polynomial type of
degree $\dim A-1$. Since $\dim A<\dim {\mathcal R}(I)=n+2$ ( see \cite[1.1(i)]{CHTV}), the
proof of the theorem is completed.
\end{proof}

We close this section with the following result whose proof follows the same line of arguments
as in the part of the proof of Theorem \ref{graded} where it is shown that $P_{I^kM}(ak)$
is of polynomial type of degree $\leq n$.

Let $A=K[x_1,\ldots, x_n, y_1,\ldots, y_m]$ be the standard
bigraded polynomial ring, that is, we have degree $x_i=(1,0)$ and
degree $y_j=(d_i,1)$ for all $i$ and $j$. Furthermore, let $E$ be a
finitely generated graded bigraded $S$-module. For all $k$ we set
$E_k=\Dirsum_iE_{i,k}$. Then each $E_k$ is a finitely generated
graded $S$-module, where $S=K[x_1,\ldots, x_n]$, and we may
consider its Hilbert polynomial $P_{E_k}(x)$ and its Hilbert
coefficients $e_i(E_k)$.

\begin{Theorem}
\label{new}  For $k\gg 0$ the dimension of $E_k$ is constant, say $d$. Then for
$i=0,\ldots, d$, the $i$th Hilbert coefficient  $e_i(E_k)$ as a
function of $k$ is of polynomial type of degree $\leq n+m-d+i$.
\end{Theorem}

\section{Asymptotic behavior of multiplicity; the graded case}
Let $J\subset I\subset S$ be graded ideals, $M$ graded $S$-modules
and $N\subset M$ a graded submodule of $M$. Then $J^kN\subset
I^kM$. Since by Theorem \ref{graded} the Hilbert coefficients of
$M/I^kM$ and $N/J^kN$ as functions of $k$ are of polynomial type,
it follows from the equation
\begin{eqnarray}
\label{pol}
 P_{I^kM/J^kN}(x)=P_{M/N}(x)+P_{N/J^kN}(x)-P_{M/I^kM}(x)
\end{eqnarray}
that the degree of $P_{I^kM/J^kN}(x)$ is constant for large $k$.
In other words, there exists an integer $k_0$ such that
$\dim(I^kM/J^kN)=d$ for all $k\geq k_0$.

More generally, if $(M_k)_{k\geq 0}$ is family of finitely generated $S$-modules with the property that there exists an integer $k_0$ such that $\dim M_k$  is constant for all $k\geq k_0$, then we call $\dim M_{k_0}$ the {\em limit dimension} of $(M_k)_{k\geq 0}$.

Let $(a_k)_{k=1,2,\ldots}$ be a sequence of real numbers. We write
$\lim_{k\to\infty}a_k\in\QQ$ to indicate that
$\lim_{k\to\infty}a_k$ exists and that the limit is a rational
number.

\begin{Proposition}
\label{limit} Let  $d$ be the limit dimension of  $(I^kM/J^kN)_{k\geq 0}$.
Then
$$\lim_{k\to \infty}\frac{e_0(I^kM/J^kN)}{k^{n-d}}\in\QQ.$$
\end{Proposition}
\begin{proof} Let $d_1=\dim M/IM$, $d_2=\dim N/JN$ and  $c=\dim M/N$.
Then $\dim M/I^kM=d_1$ and $\dim N/J^kN=d_2$ for all $k$. It
follows from equation (\ref{pol}) that
\[
e_0(I^kM/J^kN)=e_{c-d}(M/N)+e_{d_2-d}(N/J^kN)-e_{d_1-d}(M/I^kM)\quad\text{for}\quad
k\gg 0.
\]
Hence Theorem \ref{graded} implies that $e_0(I^kM/J^kN)$ is of
polynomial type of degree $\leq n-d$, and this implies that
$\lim_{k\to \infty}e_0(I^kM/J^kN)/k^{n-d}$ exists and is a
rational number.
\end{proof}

Let $\overline{J}$ denote the integral closure of an ideal $J$. As
an immediate consequence of Proposition \ref{limit} we obtain

\begin{Corollary}
\label{closure} Let $I\subset S$ be a graded ideal. Then
$(\overline{I^k}/I^k)_{k\geq 0}$ admits a limit dimension, say $d$, and
\[
\lim_{k\to \infty}\frac{e_0(\overline{I^k}/I^k)}{k^{n-d}}\in\QQ.
\]
\end{Corollary}

\begin{proof} Since $A=\Dirsum_{k\geq 0}\overline{I^k}$ is equal to the
integral closure of the Rees algebra ${\mathcal R}(I)$, it follows
that $A$ is a finitely generated graded ${\mathcal R}(I)$-module.
In particular it follows that there exists an integer $k_0$ such
that $A_k=I^{k-k_0}A_{k_0}$ for all $k\geq k_0$. In other words,
we have
$$
\overline{I^k}/I^k=I^{k-k_0}\overline{I^{k_0}}/I^{k-k_0}I^{k_0}\quad
\text{for  all}\quad  k\geq k_0.
$$
Thus the assertion follows from Proposition \ref{limit}.
\end{proof}

We give another application of Proposition \ref{limit}.

\begin{Corollary}
\label{socle}  Let $I\subset S$ be a monomial ideal and let $\mm =(x_1,\ldots,x_n)$ be the graded  maximal ideal of $S$.
Then
$((I^k\: \mm)/I^k)_{k\geq 0}$ admits a limit dimension, say $d$, and
\[
\lim_{k\to \infty}\frac{\length((I^k\:\mm)/I^k)}{k^{n-d}}\in\QQ.
\]
\end{Corollary}

\begin{proof} We will show that $M=\Dirsum_{k\geq 0} I^k\:\mm$ is a finitely generated ${\mathcal R}(I)$-module. Then we proceed as in the proof of  Corollary \ref{closure}.

We first notice that $I^k\: \mm=\Sect_{i=1}^n(I^k\: x_i)$. This implies that $M=\Sect_{i=1}^nM_i$ where $M_i=\Dirsum_{k\geq 0} I^k\: x_i$. Since ${\mathcal R}(I)$ is Noetherian it suffices  therefore to show that each $M_i$ is a finitely generated ${\mathcal R}(I)$-module. In fact, we will show that $M_i$ is generated by its degree 1 component. In other words, we will show that $I^k\: x_i=(I^{k-1}\: x_i)I$ for all $k> 1$.

Since $I$ is a monomial ideal we can write $I=Jx_i+L$ where $J$ and $L$ are again monomial ideals, and where none of the generators of $L$ is divisible by $x_i$. Then one has
\[
I^k=\sum_{j=1}^kJ^jL^{k-j}x_i^j+L^k, \quad\text{and hence}\quad I^k\: x_i=\sum_{j=1}^kJ^jL^{k-j}x_i^{j-1}+L^k.
\]
Similarly we have $I^{k-1}\: x_i=\sum_{j=1}^{k-1}J^jL^{k-1-j}x_i^{j-1}+L^{k-1}$. It follows that
\begin{eqnarray}
\label{power}
(I^{k-1}\: x_i)I&=&(\sum_{j=1}^{k-1}J^jL^{k-1-j}x_i^{j-1}+L^{k-1})(Jx_i+L)\\
&=& A+B+L^{k-1}Jx_i+L^k,\nonumber
\end{eqnarray}
with
\[
A= \sum_{j=1}^{k-1}J^{j+1}L^{k-1-j}x_i^j\quad \text{and}\quad B=\sum_{j=1}^{k-1}J^{j}L^{k-j}x_i^{j-1}.
\]
Observe that $L^{k-1}Jx_i$ appears as a multiple of a summand in $B$, so we may omit it  in the right hand sum of (\ref{power}). Next we may rewrite the sum $A$ as $\sum_{j=2}^{k}J^{j}L^{k-j}x_i^{j-1}$. Then it follows from (\ref{power})  that
\[
 (I^{k-1}\: x_i)I=\sum_{j=1}^kJ^jL^{k-j}x_i^{j-1}+L^k=I^k\: x_i,
 \]
 as desired.
\end{proof}

Let $I$ and $J$ be graded ideals in $S$. The $k$th symbolic power  of $I$ with
respect to $J$ is defined to be the graded ideal $I_k(J)=I^k:J^\infty$.

We want to single out two cases of interest: in the first case, let  $J$ be the graded maximal ideal of $S$.  Then $I_k(J)$ is the $k$th saturated power of $I$, denoted $(I^k)^{\sat}$.

In the second case, consider the set $A^*(I)=\Union_k\Ass(S/I^k)$ of asymptotic prime ideals of $I$. This set is known to be  finite (see \cite{Br}). Of course the set on minimal prime ideals $Min(I)$ of $I$ is a subset of $A^*(I)$. Let
\[
J=\Sect_{P\in\Ass^*(I)\setminus \Min(I)}P.
\]
For this choice of $J$ one obtains the ordinary symbolic powers $I^{(k)}$ of $I$.

We are interested in the limit behavior of $e_0((I_k(J)/I^k)$. In \cite[2.2]{CHST} there is given an example that shows that  $\length((I^k)^{sat}/I^k)$ need not to be of polynomial type, and not even of quasi-polynomial type. However we will see that if  $I$ and $J$ are monomial ideals, then the limit behavior of $e_0(I_k(J)/I^k)$ is quite nice. This is mainly due to the fact that in  this situation, the graded ring $\Dirsum_{k\geq 0} I_k(J)$ is a finitely generated $S$-algebra, as  shown in \cite[3.2]{HHT}.

More generally, suppose that we are given  two families $(I_k)_{k\geq 0}$ and $(J_k)_{k\geq 0}$ of graded ideals of $S$ with the property that $J_kJ_\ell\subset J_{k+\ell}$ and $I_kI_\ell\subset I_{k+\ell}$ for all $k$ and $\ell$, and such that the algebras $A=\Dirsum_{k\geq 0}J_k$ and $B=\Dirsum_{k\geq 0}I_k$ are finitely generated. We also assume that $J_k\subset I_k$ for all $k$ and that $I_0=J_0=S$.

We first show

\begin{Proposition}
\label{quasi}
The numerical function $e_0(I_k/J_k)$ is of quasi-polynomial type.
\end{Proposition}

\begin{proof}
Since $A$ and $B$ are finitely generated graded $S$-algebras, there exist integers $s$ and $t$ such that
the Veronese algebra $A^{(s)}$ of $A$ and the Veronese subalgebra $B^{(t)}$ of $B$ is standard graded, see for example \cite[2.1]{HHT}. We may assume that $s=t$. Otherwise we replace $s$ as well as $t$ by $st$.

As an $A^{(s)}$-module the algebra $A$ is a direct sum of the graded $A^{(s)}$-modules
$A^{(s;i)}=\Dirsum_{k\geq 0}A_{ks+i}$. Indeed, $A=\Dirsum_{i=0}^{s-1} A^{(s;i)}$.  Similarly, $B=\Dirsum_{i=0}^{s-1} B^{(s;i)}$.

Since $A^{(s;i)}$ is a finitely generated $A^{(s)}$-module, there exist an integer $v_i$ such that $A_{(k+v_i)s+i}=A_{ks}A_{v_is+i}$ for all $k\geq 0$. Let $v=\max\{v_i\:\; i=0,\ldots,s-1\}$. Replacing each $v_i$ by $v$ we may assume that all $v_i$ are equal to $v$. Similarly there exists an integer $w$ such that $B_{(k+w)s+i}=B_{ks}B_{ws+i}$ for all $k\geq 0$ and all $i=0,\ldots, s-1$. Taking the maximum of $v$ and $w$ we may as well assume that $v=w$.

Thus we may assume that $A$ and $B$ are standard graded and have to show that $e_0(I_1^kM/J_1^kN)$ is of polynomial type, where $M$ and $N$  are graded $S$-modules. But this has already been shown in Proposition \ref{limit}.
\end{proof}

A priori it is not clear that $\dim I_k/J_k$ is constant for $k\gg 0$. But if we assume that $A$ is standard, that is  $J_k=J^k$ for all $k$ where $J=J_1$, then we have

\begin{Theorem}
\label{tony} With the notation introduced, suppose that $A$ is standard graded. Let $P_0,\ldots, P_{g-1}$ be the polynomials such that $e_0(I_{mg+i}/J^{mg+i})=P_i(m)$ for $m\gg 0$. Then
\begin{enumerate}[\rm (a)]
\item $\dim  I_k/J^k$ is constant for $k\gg 0$.
\item   $\deg P_i=\deg P_0$ for all $i$ and all of them have the same leading coefficient.
\end{enumerate}
\end{Theorem}

\begin{Corollary}
 \label{conclusion}
 Let $I$ and $J$ be monomial ideals, and let $d$ be the limit dimension of  $(I_k(J)/I^k)_{k\geq 0}$. Then
\[
\lim_k\frac{1}{k^{n-d}} \ e_0\left(\frac{I_k(J)}{I^k}\right)\in \QQ.
\]
In particular, $$\lim_k\frac{1}{k^{n}} \ \ell\left(\frac{(I^k)^{\sat}}{I^k}\right)\in \QQ.$$
\end{Corollary}

\medskip
In order to prove Theorem \ref{tony} we  need to introduce some definitions.
Let $M=\Dirsum_{k\geq 0}M_k$ be a  bigraded $A=\Dirsum_kJ^k$-module over the bigraded algebra $A$ such that each $M_k$ is a finitely generated $S$-module.  All modules considered in the sequel will be of this kind.

The module $M$ is called {\em quasi-finite}, if $H^0_{A_+}(M)_i=0$ for $i\gg 0$. Here  $H^0_{A_+}(M)$ denotes the $0$th local cohomology of $M$ with respect to the ideal $A_+=\Dirsum_{k>0}A_k$. We consider quasi-finite modules because in the application we have in mind the modules $\Dirsum_{k\geq 0}I_k(J)/I^k$ are in general not finitely generated $\Dirsum_kI_k$ modules, but only quasi-finite, as we shall see.

We say an   element $x \in A_1 $ is {\em $M$-filter regular} if  $(0 \colon_M x)_i = 0$ for $i \gg 0$. We
do not require  that $x$ is bihomogeneous.
The existence of an $M$-filter regular element is useful  due to the following:
\begin{Proposition}\label{equaldegree}
Let $M = \Dirsum_{k \geq 0}M_k$ be an $A$-module such that $e_0(M_k)$  is  a quasi-polynomial. Let $P_0,\ldots, P_{g-1}$ be the polynomials such that $e_0(M_{mg+i})=P_i(m)$ for $m\gg 0$.
If there exists an element in $A_1$ which is $M$-filter regular, then
\begin{enumerate}[\rm (a)]
\item $\dim M_k$ is constant for $k\gg 0$.
\item $\deg P_i=\deg P_0$ for all $i$ and all of them have the same leading coefficient.
\end{enumerate}
\end{Proposition}

\begin{proof}
(a) Suppose for
$m \geq s$ we have
$e_0(M_{mg+i}) = P_i(m)$ for $i =0,\ldots g-1$. We choose  $x \in A_1$  which is  $M$-filter regular. Then the  map
$M_{j}\to  M_{j+1}$ induced by multiplication with $x$ is injective for all
$j \gg 0$.
This implies that for $m \gg 0 $ we have
\[
\dim M_{mg} \leq \dim  M_{mg+1} \leq \ldots \leq  \dim M_{mg + g-1} \leq
\dim M_{(m+1)g}.
\]
The assertion follows.

(b) We assert that $e_0(M_k)\leq e_0(M_{k+1})$ for all $k\gg 0$. Let $d$ be the limit dimension of the family $(M_k)_{k\geq 0}$. If $d=0$, then $e_0(M_k)=\ell(M_k)$ and the assertion is obvious, because the  maps $M_k \to M_{k+1}$ induced by multiplication by $x$ are injective for $k\gg 0$.

So now let $d>0$, and let $c$ be the initial degree of the element $x$ in the graded module $A_1$. Then the maps
$M_k \to M_{k+1}$ induce maps $M_{k, i\geq s} \to M_{k+1, i\geq s+c}$ for all $s$, where $M_{k,i\geq t}=\Dirsum_{i\geq t}M_{k,i}$.

Notice that $e_0(M_k)=e_0(M_{k, i\geq s}) \forall s$  and that both modules, $M_{k, i\geq s}$ and $M_{k+1, i\geq s+c}$ are generated in a single degree for $s$ large enough. We choose such an $s$. Let $\mm$ be the graded maximal ideal of $S$. Since $M_{k, i\geq s}$ is generated in a single degree, the associated graded module $\gr_\mm\big((M_{k, i\geq s})_\mm\big)$ is isomorphic to $M_{k, i\geq s}$, up to a shift. Hence $e_0(\mm, (M_{k, i\geq s})_\mm)=e_0(M_{k, i\geq s})$. On the other hand we have that $(M_{k, i\geq s})_\mm\to (M_{k+1, i\geq s+c})_\mm$ is injective for all $k\gg 0$. Since both modules have the same dimension, the inequality in the next display follows from  \cite[4.7.7]{BH}. $$e_0(M_k)=e_0(\mm, (M_{k, i\geq s})_\mm)\leq e_0(\mm, (M_{k+1, i\geq s+c})_\mm)=e_0(M_{k+1}).$$

 Thus for $m \geq s$ we have that
\begin{equation}\label{ineq}
P_0(m) \leq P_1(m)\leq \cdots \leq P_i(m)\leq \cdots \leq  P_{g-1}(m) \leq P_0\left(m+1\right).
\end{equation}
The inequalities in (\ref{ineq}) imply that all $P_i$ have the same degree.
Let $s_i$ be the leading coefficient of $P_i$.  As all $P_i$ have the same degree note that
(\ref{ineq}) implies that
\[
s_0  \leq  s_1   \leq  \cdots  \leq  s_{g-1}   \leq  s_0.
\]
Thus $s_0 = s_1 =  \cdots = s_{g-1}$.
\end{proof}

Next we discuss the existence of a filter regular element for a quasi-finite module.
Set $N = M/H^{0}_{A_+}(M)$. If $M$ is finitely generated as an $A$-module,  then $x\in A_1$ is $M$-filter regular precisely when $x$ is $N$-regular, see \cite[2.1]{t}. Thus $M$-filter regular elements always  exists if $K$ is infinite.    For quasi-finite modules we have
\begin{Lemma}
\label{suplemma}
Let  $M = \bigoplus_{k \geq 0}M_k$ be a quasi-finite $A$-module.
 Set $N = M/H^{0}_{A_+}(M)$. Then
 \begin{enumerate}[\rm (a)]
 \item
 $\Ass_{A} N \cap \Var(A_+) = \emptyset$.
 \item
 $\Ass_{A} N$ is a countable set.
 \item
 If $x \in A_1$ is $N$-regular, then it is $M$-filter regular.
 \item
 If $K$ is uncountable, then there exists $x \in A_1$ which is
 $N$ regular and so $M$-filter regular.
 \end{enumerate}
\end{Lemma}
\begin{proof}
 (a) Suppose there exists $\pp \in \Ass_{A} N \cap \Var(A_+)$. Let $\pp = (0 \colon \overline{x})$ where $\overline{x} = x + H^{0}_{A_+}(M)$ is a nonzero homogeneous
element of $N$. As $\pp \overline{x} = 0$ we have $E = \pp x \subseteq H^{0}_{A_+}(M)$. Notice $E$ is finitely generated as an $A$-module. So $A_{+}^{m}E = A_{+}^{m}\pp x =0$. Since $ \pp \supseteq A_{+}$ we have that
$A_{+}^{m+1}x =0$. Therefore $x \in H^{0}_{A_+}(M)$, a contradiction.

(b) For $i \geq 0$, set
$D_i$ to be the $A$-submodule of $N$ generated by $N_0,\ldots, N_i$. As each $N_i$ is a finitely
generated $S$-module, it follows that each $D_i$ is a finitely generated $A$-module. The assertion follows since
$N = \bigcup_{i\geq 0} D_i$.

(c) Let $x$ be $N$-regular.  As $M$ is quasi-finite we have $M_i = N_i$ for $i \gg 0$. It follows that
$(0 \colon_M x)_i = 0$ for $i \gg 0$. Thus $x$ is $M$-filter regular.

(d)  By (b)  we get that $\Ass_{A_+} N$ is countable. Say $\Ass_{A} N = \{ Q_1,Q_2,\ldots \}$. Then for each $i \geq 1$ we have that $Q_i \cap A_1$ is a
proper graded $S$-submodule of $A_1$. Since $K$ is uncountable,  there exists $x \in A_1 \setminus \bigcup_{i \geq 1} Q_i$. Clearly
$x$ is $N$-regular.
\end{proof}

\begin{proof}[Proof  of \ref{tony}]
We may assume $K$ is uncountable. If not,  we choose a field extension $K^{\prime}$ of $K$ which is uncountable, and  set $S^\prime  = S \otimes_{K} K^\prime$, $A^\prime = A\otimes_{S} S^\prime$ and $B^\prime = B\otimes_{S} S^\prime$.  It can be checked that $e_0(B^{\prime}_k/
 A^{\prime}_{k}) = e_0(B_k/ A_k)$ for all $k \geq 1$.

 Set $M = B/A$.
 The exact sequence $0 \rightarrow A \rightarrow B \rightarrow B/A \rightarrow 0$ yields an exact sequence
\[
0 \rightarrow H^{0}_{A_+}(A)  \rightarrow H^{0}_{A_+}(B)  \rightarrow H^{0}_{A_+}(B/A) \rightarrow H^{1}_{A_+}(A) \rightarrow
\]
We assert that $B/A$ is a quasi-finite $A$-module. In fact, $H^0_{A_+}(B)= H^0_{A_+B}(B)=0$ because $B$ is a domain, and $H^{1}_{A_+}(A)_k = 0$ for $k\gg 0$ in general \cite[15.1.5]{bs}. Thus the long exact cohomology sequence implies that $H^{0}_{A_+}(B/A)_k=0$ for $k\gg 0$, as desired.
  By Lemma \ref{suplemma}(d) we have
 that there exists $M$-filter regular elements. So by Proposition \ref{equaldegree} the result follows.
  \end{proof}

\section{Asymptotic behavior of multiplicity in the local case}
 Let $(R,\mm)$ be a local ring. If $M$ is a finitely generated $R$-module then
 $e_0(M)$ denotes the Hilbert-Samuel multiplicity of $M$ with respect to $\mm$. Let  $\mathcal{F} = \{I_n \}_{n \geq 0} $ be a Noetherian filtration on $R$ i.e., a chain of ideals
$$R = I_0 \supseteq I_1\supseteq \cdots \quad \supseteq I_i \supseteq I_{i+1} \supseteq \cdots \text{such that } \ I_iI_j  \subseteq I_{i+j} \  \  \text{for all $i$,$j$ and } $$
 $\Rees(\Fc) = \bigoplus_{k \geq 0} I_k$ is a graded (not necessarily standard) Noetherian $R$-algebra.
\begin{Example}
 \emph{Let $I^{(k)}$ denote  the $k^{th}$ symbolic power of $I_1$. There are many cases known when $\mathcal{F} = \{I^{(k)} \}_{k \geq 0} $
is Noetherian ;  see for example \cite{C} \cite{H} and \cite{GNW}.  The symbolic power filtration need not be Noetherian  in general \cite{rob}.}
\end{Example}

The main result of this section is the local analogue of Theorem \ref{tony}
\begin{Theorem}\label{filt}
Let $J$ be an ideal  in $R$ with positive grade such that $ J \subseteq I_1$. Then
\begin{enumerate}[\rm (a)]
\item
$\dim I_k/J^k$ is constant for $k \gg 0$.
\item
There exists polynomials $P_0, \ldots, P_{g-1}$ such that $P_i(m) = e_0(I_{mg + i}/J^{mg +i})$ for all
$m \gg 0$. In other words,  $e_0(I_{k}/J^{k})$ is of quasi-polynomial type  of period $g$.
\item
$\deg P_i = \deg P_0$ for all $i$ and all of them have the same leading coefficient.
\end{enumerate}
\end{Theorem}
It is convenient to prove  Theorem \ref{filt} by considering the following more
general setup: let $A \subseteq B$ be a homogeneous  inclusion of Noetherian graded algebras with
$R = A_0 = B_0$ a local ring with maximal ideal $\mm$. We assume $A$ is a standard graded $R$-algebra,  while $B$ is {\em not necessarily} a standard graded $R$-algebra.

A typical example of this situation is when $B=\Rees(\Fc)$ and $A=\Rees(J)$ with $J\subset I_1$.

In analogy to the graded case we call a graded $A$-module $M=\Dirsum_{k\geq 0}M_k$ {\em quasi-finite} if
each $M_k$ is a finitely generated $R$-module and $H^{0}_{A_+}(M)_j = 0$ for all $j\gg 0$.

Note that $B$ is  trivially a  quasi-finite $A$-module if $\grade(A_1B, B) > 0$. Notice if  $A_n, B_n$ are ideals in $R$ with  $\grade(A_1,R)> 0$, then $\grade(A_1B, B) > 0$, and so $B$ is  a  quasi-finite $A$-module.

In this more general frame Theorem \ref{filt} reads as follows:

\begin{Theorem}\label{MTinclusion}
Assume that $B$ is quasi-finite $A$-module. Then
\begin{enumerate}[\rm (a)]
\item
$\dim B_k/A_k$ is constant for $k \gg 0$.
\item
There exists polynomials $P_0, \ldots, P_{g-1}$ such that $P_i(m) = e_0(B_{mg + i}/A_{mg +i})$ for all
$m \gg 0$. In other words,  $e_0(B_{k}/A_{k})$ is  of quasi-polynomial type  of period $g$.
\item
$\deg P_i = \deg P_0\leq \dim B+d-1$ for all $i$ and all of them have the same leading coefficient. Here $d$ is the limit dimension of $(B_k/A_k)_{k\geq 0}$.
\end{enumerate}

\end{Theorem}

We give an example which shows that the conclusion of the theorem would fail without the quasi-finite condition on $B$.
\begin{Example}
{\em Let $R=K[[x]]$ be the power series ring over a field $K$. Set   $A = R[u,v]$  and $B = R[u,v,w]/(wu,wv)$. We give $w$ degree 2.
Thus $wB \subseteq H^{0}_{A_+}(B )$. Thus $B$ is not a quasi-finite $A$-module.
It can be
easily seen that $B_k/A_{k} = 0$ if $k$ is odd and $Rw^{k/2}$ if $k$ is even. So $g = 2$ and the multiplicity polynomial is
$(0,1)$.}
\end{Example}

One of the chief problems in trying to prove the above theorem is that $B/A$ is not a $B$-module. It is also, in general,
 not a
finitely generated $A$-module.
Rees's technique \cite{Rees1}, in this case   is to consider the Rees ring $S = \Dirsum_{k\geq 0} (A_+B)^k$ and the
associated graded ring $G = \gr_{A_1B} B $ of the  ideal $A_1B$ in $B$.

Observe that $(A_+B)^k=\Dirsum_{i\geq 0} A_kB_i$, so that $S=\Dirsum_{i,j}A_iB_{j-i}$. Hence $S$ is a naturally bigraded $R$-algebra with $S_{ij}=A_iB_{j-i}$ for all $i,j$. We have $S_{i,j} = 0$ if $j<i$ and $S_{i,i} = A_i$.

Notice that $G$ inherits a bigraded structure from $S$, so that $G = \bigoplus_{i \geq 0,j\geq 0}G_{i,j}$ with
$$G_{i,j} = \frac{S_{i,j}}{S_{i,j+1}} = \frac{B_{i-j}A_j}{B_{i-j-1}A_{j+1}}.$$
$G$ also has the property that $G_{i,j} = 0$ if $j<i$ and $G_{i,i} = A_i$. If we consider $A$ as bigraded
algebra with $A_{i,j} = 0$ if $i \neq j$ and $A_{i,i} = A_i$ then $A$ becomes a bigraded $R$-subalgebra of $G$.
Finally let $L = \bigoplus_{j> i}G_{i,j}$. Notice that
$L$ is an ideal of $G$.

We now consider $G$, $L$ singly graded by summing across columns, i.e, set $G = \bigoplus_{i\geq 0}G_i$, $L =\bigoplus_{i\geq 0}L_i $ where
$$     G_{i} = \bigoplus_{j \geq 0} G_{i,j} =  \bigoplus_{j = 0}^{i} G_{i,j} \quad \text{and} \quad     L_{i} = \bigoplus_{j \geq 0} L_{i,j}.     $$

The advantage of  Rees's technique is due to the following observation:
consider the filtration of   $R$-modules
\begin{eqnarray}
\label{f1}
 A_i \subseteq B_1 A_{i-1} \subseteq B_2 A_{i-2}\subseteq B_{i-1} A_{1} \subseteq B_{i}.
 \end{eqnarray}
One can check that
\begin{eqnarray}
\label{f2}
L_i = \bigoplus_{j = 0}^{i-1} \frac{B_{i-j} A_{j}} {B_{i-(j+1)} A_{j+1}}.
\end{eqnarray}

This leads to the following

\begin{Lemma}\label{advtReesTech}
With the notation introduced we have an inclusion of graded algebras $A \subseteq G$ and a finitely generated $G$-module $L$ with
\begin{enumerate}[\rm (a)]
\item
$\dim_R B_k/A_k = \dim_R L_k$ for all $k \geq 1$.
\item
$e_0(B_k/A_k) = e_0(L_k)$ for all $k \geq 1$.
\item
$\dim (G) = \dim (B)$.
\item
$e_0(B_k/A_k)=e_0(L_k)$ is of quasi-polynomial type.
\end{enumerate}
\end{Lemma}

\begin{proof}
The first two statement follow from (\ref{f1}) and (\ref{f2}). For the third statement see \cite[4.5.6]{BH}
Finally, statement (d) follows from Proposition~\ref{XX}.
\end{proof}

Note that assertion (b) of Theorem \ref{MTinclusion} follows from statement (d) above. The remaining assertions of  Theorem \ref{MTinclusion} depend  on the existence of a $B/A$-filter regular element in $A_1$. One proves its existence just as in the graded case by using a local version of Proposition~\ref{equaldegree} and Proposition~\ref{suplemma}.

\begin{proof}[Proof  of \ref{MTinclusion} (a) and (c)]
If $K = R/\mm$ is not uncountable then we do the following standard trick. First we complete $R$. Then we consider
$R^\prime = \widehat{R}[[x]]_{\mm[[x]]}$ where $\widehat{R}$ is the completion of $R$. Notice that $R^\prime$ has residue field $K((x))$, which is uncountable.
Set $A^\prime = A\otimes_{R} R^\prime$ and
$B^\prime = B\otimes_{R} R^\prime$. Since $A^\prime$ is a flat extension of $A$ it follows that
  $B^\prime$ is also a quasi-finite $A^\prime$-module. It can be checked that $e_0(B^{\prime}_k/
 A^{\prime}_{k}) = e_0(B_k/ A_k)$ for all $k \geq 1$.

Similarly as in the proof of Theorem~\ref{tony} it follows that  $B/A$ is a quasi-finite $A$-module.
 By the local version of Proposition~\ref{suplemma}(d) there exists a $B/A$-filter regular element in $A_1$ . So by the local version of Proposition~\ref{equaldegree} the results follow, except the degree bound for the $P_i$. But this is a consequence of Proposition~\ref{XX}.
\end{proof}

  Our main theorem is an easy corollary of  Theorem \ref{MTinclusion}.

\begin{proof}[Proof of Theorem \ref{filt}]
Let $A = \Rees(J)$ and $B = \Rees(\Fc)$. Since $J \subseteq I_1$ we have an
inclusion of graded algebras $A \subseteq B$. Since $\grade_R J>0$ we get that
 $\grade_BA_1B>0$. So $H^{0}_{A_+}(B)=H^{0}_{A_1B}(B) = 0$.
In particular,  $B$ is a quasi-finite $A$-module. The result
now follows from Theorem \ref{MTinclusion}.\end{proof}

\section{The Rees polynomial of $I$ and $J$}
Throughout this section, let $(R, \mm)$ be a local ring of dimension $d.$ For an $\mm$-primary 
ideal $I,$ put $H(I,n)=\ell(R/I^n)$  where $\ell$ denotes length. Let $P(I,n)$ denote the Hilbert 
polynomial corresponding to the Hilbert function $H(I,n).$ We write
$P(I,n)$ in the form:
$$P(I,n)=e_0(I) \binom{n+d-1}{d}-e_1(I) \binom{n+d-2}{d-1}+\cdots+(-1)^de_d(I).
$$
Recall that an ideal $J\subset I$ is called a reduction of $I$
if  $JI^n=I^{n+1}$
for some positive integer $n.$ The Rees multiplicity theorem asserts that
if $R$ is a  quasi-unmixed local ring then $J$ is a reduction of $I$ if and only if 
$e_0(I)=e_0(J).$ The Rees multiplicity theorem has been generalized
for ideals that are not $\mm$-primary by a number of authors. In \cite{r}
Rees provided one such generalization which we now state.
Let $J \subset I$ be ideals of $R$ such that $\ell(I/J)$ is finite.
Then the numerical function
$H(I/J,n)=\ell(I^n/J^n)$ is  given by a polynomial  $P(I/J,n)$ for large $n,$ \cite{a}.
We call $H(I/J,n)$        and $P(I/J,n)$ the {\it Rees function} and the
{\it Rees polynomial}  of the pair
$(J,I)$ respectively.
Rees showed that if $(R,\mm)$ is quasi-unmixed then $J$  is a reduction of $I$
if and only if $\deg P(I/J,n) < d.$  It is natural to ask for the exact degree
of $P(I/J,n).$ In case $I$ and $J$ are $\mm$-primary it is easy to see
that $\deg P(I/J,n)< d-k $ if and only if  $e_i(I)=e_i(J)$ for $i=0,1, \ldots,k$ for all $k=0,1, 
\ldots, d.$  K. Shah  \cite{s} introduced coefficient ideals
which characterize  the degree of $P(I/J,n)$ when $I$ and $J$ are $\mm$-primary.

\begin{Theorem}[K. Shah] Let $(R,\mm)$ be a quasi-unmixed local ring of positive dimension
$d$ with infinite residue field $R/\mm.$ Let $I$ be an $\mm$-primary ideal.
Then there exist unique largest ideals $I_1, I_2, \ldots, I_d$ containing $I$
such that $e_i(I)=e_i(I_k)$ for $i=0,1, \ldots,k$ and
$I \subset I_d \subset I_{d-1} \subset \cdots I_1 \subset I_0=\bar{I}. $ The ideal
$I_k$ is  called the $kth$-coefficient ideal of $I.$
\end{Theorem}

We will  prove an analogue of the above theorem for non-$\mm$-primary ideals
which will help us in determining the degree of $P(I/J,n).$ The proof
is similar to the proof of Shah's theorem. First we need the following

\begin{Definition}  Let $(R,\mm)$ be a local ring and let $I$ be an ideal in $R$.
Let $I^{sat} = \bigcup_{n\geq 1}(I \colon\mm^n)$ be the saturation of $I$.
We call $\qq(I) = \overline{I}\cap I^{sat}$ the relative
integral closure  of $I.$
\end{Definition}

In several results in this section  we consider ideals $J \subseteq I$ so that
$J$ is a reduction of $I$ and $\ell(I/J) < \infty.$  For these cases we have:

\begin{Proposition}
Let $(R,\mm)$ be a local ring and $I$ an ideal of $R.$ Then

\begin{enumerate}
\item
$I$ is a reduction of $\qq(I)$ and $\ell(\qq(I) /I) < \infty$.
\item
If $I$ is a reduction of $K$ with $\ell(K /I) < \infty$ then $K \subseteq \qq(I).$
\item
If $J \subseteq I$ then $\qq(J)\subseteq  \qq(I)$.
\item
$\qq(\qq(I) ) = \qq(I)$.
\end{enumerate}
\end{Proposition}

\begin{proof}
The assertions (1), (2), (3) are easy to show. We prove (4). If $\qq(I)$ is a reduction of $K$ with $\ell(K /\qq(I)) < \infty$ then note that $K$ is also a reduction of $I$ and also $\ell(K /I) < \infty$. So by (2) we have
$K \subseteq \qq(I).$ Thus $K = I$. By taking $K = \qq(\qq(I) )$ we get the desired result.
\end{proof}
We now give an analogue of Shah's Theorem.
\begin{Theorem} \label{coefficient}
 Let $(R,\mm)$ be a $d$-dimensional quasi-unmixed local ring
with  infinite residue  field $R/\mm.$ Let $I$ be any ideal of analytic spread $a.$
Then for $k=1, 2, \ldots, a,$
there  exist unique largest ideals $I_k$  with  $\ell(I_k/I) < \infty, $
and $I \subseteq I_a \subseteq I_{a-1} \subseteq  \cdots \subseteq I_1 \subseteq \qq(I), $
such that $\deg P(I_k/I,n)< a-k $  for $k=1, 2, \ldots, a.$
\end{Theorem}

\begin{proof} For  each $k=1, 2, \ldots, a,$ consider the sets
$$ V_k=\left\{ L \mid L \; \mbox{is an ideal of $R$ with }
              \ell(L/I) < \infty,\;  \deg P(L/I,n) < a-k \right\}.$$
If $L \in V_k,$ then $ \deg P(L/I,n) < d-1.$ Hence by Rees's theorem,
$I$ is a reduction of $L.$ Hence $L \subseteq \qq(I).$

Since  $I \in V_k$ and   $R$  is Noetherian, $V_k$  has a maximal member, say
$J.$ We show that $J$ is unique. Let $L \in V_k$ and $x \in L.$
  Since $I \subset (I,x) \subset L,$ for all $n \geq 1,$
$\ell( (I,x)^n/ I^n) \leq \ell(L^n/I^n) < \infty.$ Hence
$\deg P((I,x)/I,n) \leq \deg P(L/I,n) < a-k.$
By Rees's theorem, $I$ is a   reduction of $(I,x).$ Let $I(I,x)^t= (I,x)^{t+1}$
for some integer $t.$ Then $x^{t+1} \in I(I,x)^t \subset J(J,x)^t.$ Hence
$J(J,x)^t=(J,x)^{t+1}.$  It follows that  for all
$n \geq t, J^{n-t}(J,x)^t=(J,x)^n.$ Since $\ell((I,x)/I) < \infty,$ there
exists an $r$ such that $x \mm ^r \subset I \subset J.$ Hence
$\ell((J,x)/I)< \infty.$  Now we  show that $\deg P((J,x)/I, n)< a-k .$
For all $n \geq t,$ we have:

\begin{eqnarray*}
\ell((J,x)^n/I^n)&=& \ell(J^{n-t}(J,x)^t/I^n)
                 = \ell((J^n,J^{n-1}x, \ldots, x^tJ^{n-t})/I^n) \\
				 & \leq & \sum_{i=1}^t\ell \left((J^n+x^iJ^{n-i})/I^n \right)
				  \leq   \sum_{i=1}^t \left[\ell(J^n/I^n)+
				               \ell((I^n+x^iJ^{n-i})/I^n)  \right] \\
				&=& 	\sum_{i=1}^t \left[\ell(J^n/I^n)+ \ell\left((x^iI^{n-i}+I^n)/I^n \right)
				                       +\ell\left( \frac{x^iJ^{n-i}+I^n}{x^iI^{n-i}+I^n}\right) 
\right] \\
				&\leq & \sum_{i=1}^t\ell(J^{n-i}/I^{n-i})+\ell((I,x)^n/I^n)+\ell(J^n/I^n). \\          		
 \end{eqnarray*}
 Hence $\deg P((J,x)/I, n) < d-k.$ Therefore $(J,x) \in V_k.$ By maximality of $J,$ we conclude 
that $x \in J.$ Hence $L \subset J$ and therefore $J$ is  unique.
\end{proof}

\begin{Corollary} Let $(R,\mm)$ be a quasi-unmixed local ring. Let $I \subset J$ be ideals of
$R$ such that $\ell(J/I) < \infty.$ Then $\deg P(J/I,n)=a-k \; \text{if and only if}\; J \subset 
I_{k-1} \;\text{but}\; J  \not  \subset I_k $ for $k=1, 2, \ldots,a.$
\end{Corollary}

As a consequence of the next result, we   identify the degree of
the Rees polynomial of $J$ and $I.$

\begin{Proposition} Let $R=\bigoplus_{i=0}^\infty R_n$ be a standard graded Noetherian ring over a
Noetherian local ring $(R_0,\mm).$ Let $M=\bigoplus_{i=0}^\infty M_n$ be a finitely generated graded
$R$-module with $\ell(M_n) < \infty $ for all $n \geq 0.$ Let $H(M,n)=\ell_{R_0}(M_n)$ be the 
Hilbert function of $M.$   Then  $\dim M=\dim M/\mm M$ and  $H(M,n)$ is a polynomial function of 
degree $\dim M-1.$
\end{Proposition}

\begin{proof}
Since $M$ is a finitely generated $R$-module and $\ell(M_n)  < \infty,$ for all
$n \geq 0,$ there exists an $r$ such that $\mm^rM_n=0$ for all $n \geq 0.$ Thus
$\mm^r R \subset \text{ann}_R(M).$ Since $\Supp M/\mm M= \Supp M \cap V(\mm R),$ it follows that
$\Supp M =\Supp M/\mm M$ and hence $\dim M=\dim M/\mm M.$

\noindent
Let $q= \oplus_{n=0}^{\infty} q_n= \text{ann}_R (M).$ Then each $M_n$ is an $R_0/q_0$-module.
Let $M= \sum_{i=1}^g Rm_i$ where  $\deg m_i=d_i$ for $i=1, 2, \ldots, g.$ Define the $R$-module
homomorphism $\phi : R \To M^{\oplus g}$ by $\phi(r)=(rm_1, rm_2, \ldots, rm_g).$ Then
$\Ker \phi =\text{ann}M.$ Restricting $\phi$ to  $R_0,$ we see that $R_0/q_0$ is a submodule
of an $R_0$-module of finite length. Hence $q_0$ is $\mm$-primary.
Let $S=R/\Ann M=\oplus_{n=0}^{\infty}S_n.$ Then $M$ is a finitely generated $S$-module where
$S_0$ is an Artin local ring. Hence by Hilbert-Serre theorem, $H(M,n)$ is a polynomial function
of degree $\dim M-1.$
\end{proof}

\begin{Corollary}
Let $(R,\mm)$ be a local ring of dimension $d.$ Let $J \subset I$
be ideals of $R$ such that  $J$  is a reduction of $I$ and
$\ell(I/J) < \infty.$ Let  $M=\Rees(I)/\Rees(J)=\bigoplus_{n=0}^{\infty}I^n/J^n.$
Then $H(I/J,n)$ is a polynomial function of degree   $\dim M-1.$
\end{Corollary}

\begin{proof}  Since $J$ is a reduction of $I,$ $\Rees(I)$ and  $M$  are  finite
$\Rees(J)$-module. Now apply the above proposition to $M.$
\end{proof}

The maximum degree of the Rees polynomial $P(I/J,n)$ of a reduction $J$
of an ideal $I$ is $d-1.$ In the next theorem we give a sufficient condition
so that this degree is $d-1.$

\begin{Theorem} \label{S2} Let $(R, \mm)$ be a $d$-dimensional quasi-unmixed local domain.
Let $J$ be a proper reduction  of an ideal $I$ of $R$ with $\ell(I/J) < \infty.$
If the Rees ring $\Rees(J)$ satisfies Serre's condition $S_2,$ then $a(I)=d$ and
$\deg P(I/J,n)=d-1.$
\end{Theorem}

\begin{proof} Let $A=\Rees(J)$ and $B=\Rees(I).$ Then $A$ and $B$  have same
quotient field and $B$ is a finite integral extension of $A.$ Since $A$
satisfies the $S_2$ condition, the conductor $A:B$ is a height one unmixed
ideal. Indeed, let $B=Ab_1+Ab_2+\cdots+Ab_n$ for some
$b_j \in B$ for $j=1, 2, \ldots,n.$ Let $b_j=x_j/y_j$ for some $x_j, y_j \in A,$
for $j=1, 2, \ldots, n.$ Then $A:B=\cap_{j=1}^n A:b_j=\cap_{j=1}^n (y_j:x_j).$
Since $A$ satisfies the $S_2$ condition, principal ideals in $A$ have no embedded
primary components. Hence $A:B$ is a height one unmixed ideal. Since
$\Ann_A M =A:B,$ and $A$ is universally catenary,
we have $\dim M=\dim \Rees(J)/\Ann M={d+1}-1=d. $ Hence
$\deg P(I/J,n)=\dim M-1=\dim M/\mm M-1=d-1 \leq a(J)-1 \leq d-1.$ Therefore $a(J)=a(I)=d.$
\end{proof}

\begin{Remark}
The relative integral closure can equal $I$ if either $\overline{I} = I$ or $I^{sat} = I$. We call these as trivial cases.
 The following  is a nice consequence to Theorem \ref{S2}
\end{Remark}

\begin{Corollary}
Let $(R, \mm)$ be a $d$-dimensional quasi-unmixed local domain and let $I$ be an ideal in $R$. If $a(I) < d$ and
 $\Rees(I)$ is $S_2$ then $\qq(I) = I$.
\end{Corollary}
\begin{proof}
 If $\qq(I) \neq I$ then note that $I$ is a proper reduction of $\qq(I)$ and
$\ell(\qq(I)/I)$ is finite. So by Theorem \ref{S2} we get $a(\qq(I)) = a(I) = d$ which is a contradiction.
\end{proof}

The degree of the Rees polynomial of a pair of  normal ideals can be completely determined in analytically unramified quasi-unmixed local rings. This follows from the next theorem.

\begin{Theorem} \label{normalRees}
Let $(R,\mm)$ be an analytically unramified quasi-unmixed local ring of dimension $d.$
Let $J \subset I$ be ideals of $R$ such that $J$ is not a reduction of $I$ and
$\ell (I/J) < \infty.$ Let $f(n)=\ell(\ol{I^n}/\ol{J^n}).$ Then $f(n)$ is a
polynomial function of degree $d.$

\end{Theorem}

\begin{proof} Since $\ell(I/J) <  \infty,$ there is a positive integer $r$
such that $\mm^r I \subseteq J.$ Hence $\mm^{r} \ol{I} \subseteq \ol{J}.$ Therefore
$\ell(\ol{I}/\ol{J})<  \infty.$ Since $\ell(I^n/J^n) < \infty$ for all $n,$ it follows that
$f(n) < \infty $ for all $n.$ Consider the filtration of ideals
$$ \ol{J^n} \subseteq \ol{IJ^{n-1}} \subseteq \ol{I^2J^{n-2}}
    \subseteq \cdots \ol{I^{n-1}J} \subseteq \ol{I^n}.$$
Then $f(n)=\sum_{i=0}^{n-1} \ell \left(\ol{J^iI^{n-i}}/ \ol{J^{i+1}I^{n-i-1}}\right).$
The integral closure of the
bigraded Rees algebra $\Rees:=\Rees(I,J)=\bigoplus_{r,s \geq 0} I^rJ^su^rv^s$ in the polynomial ring
$R[u,v]$ is the bigraded ring  $\ol{\Rees}:=\bigoplus_{r,s \geq 0} \ol{I^rJ^s}u^rv^s.$

Since $R$ is analytically unramified, $\ol{\Rees}$ is a finite $\Rees$-module. Indeed, by Theorem 1.4 of \cite{rees2}, there exists an integer $t$ such that $\ol{I^rJ^s} \subset I^{r-t}J^{s-t}$
for all $r,s \geq t.$ For  positive integers $r$ and  $s,$ let   $V_s=\bigoplus_{r=0}^{\infty} \ol{I^rJ^s} \subset \bigoplus\ol{I^r}$ and
 $H_r=\bigoplus_{s=0}^{\infty} \ol{I^rJ^s} \subset \bigoplus \ol{J^s}.$ Then $V_s$ is finite $\Rees(I)$-module and $H_r$ is a
finite $\Rees(J)$-module by Corollary 9.2.1 of  \cite{HS}. Thus $\ol{\Rees}$ is a finite $\Rees$-module.

Now  consider the ideals   $A=\bigoplus_{r \geq 1, s \geq 0}\ol{I^rJ^s}$ and
$B= \bigoplus_{r \geq 1, s \geq 0}\ol{I^{r-1}J^{s+1}}$ of $\ol{\Rees}.$
Then $C:=A/B$ is a finitely generated  bigraded $\Rees$-module. Since
each bigraded component of $C$ has finite length, it follows that
$H(C; x,y):=\sum_{r \geq 1, s \geq 0}\ell(C_{rs})x^ry^s$ is a rational function
of the form $H(C;x,y)=h(x,y)/(1-x)^p(1-y)^q$ for some polynomial $h(x,y) \in \ZZ[x,y].$
Note that $f(n)=\sum_{i=0}^{n-1} \ell(C_{(n-i, i)}).$ Hence
\begin{align*}
\sum_{n=0}^{\infty} f(n)z^n &= \sum_{n=0}^{\infty} \sum_{i=0}^{n-1}\ell(C_{(n-i,i)})z^n \\
                            &=  h(z,z)/(1-z)^{p+q}.
\end{align*}
Hence $f(n)$ is a polynomial function. It remains to show that its degree is $d.$
Since $R$ is analytically unramified,there exists an integer $l$ such that
$(\ol{I^l})^n=\ol{I^{ln}}$ and  $(\ol{J^l})^n=\ol{J^{ln}}$ for all $n \geq 1,$
by Corollary 9.2.1 of \cite{HS}. Since $J$ is not a reduction of $I,$ $\ol{J^l}$
is not a reduction of   $\ol{I^l}.$  Hence $\ell\left( (\ol{I^l})^n/ (\ol{J^l})^n\right) =f(nl)$  is
a polynomial function in $n$ of degree $d$ by   \cite{r}.
Hence $f(n)$ is a polynomial function of  degree $d.$
\end{proof}

The
following example shows that we cannot hope for a \emph{relative minimal reduction} in general.
\begin{Example}
Let $(R, \mm)$ be a $d$-dimensional Cohen-Macaulay ring and infinite residue field. Let $I$ be an ideal with
 a minimal reduction $W = (u_1,\ldots, u_s) $ which is generated by a regular sequence. Assume $0< s< d$.
 Then if $J$ is a reduction of $I$ with $\ell (I/J)$ finite then there exists a reduction $K$ of $I$ with
 $K \subseteq J$, $K \neq J$ and $\ell(I/K) < \infty$.
\end{Example}
\begin{proof}
Let $L$ be a minimal reduction (in the usual sense) of $I$ such that $L \subseteq J$. Note that $L$ is also generated by a regular sequence. Notice as $R/L$ is Cohen-Macaulay of positive dimension we get in particular
$L^{sat} = L$. In particular $\qq(L) = L$ and so $\ell(I/L) = \infty$.

As $\ell (I/J)$ is finite there exists $c$ such that $\mm^cI \subseteq J$. Consider the sequence of ideals
$J_n = L + \mm^{n}I$ where $n \geq c$. Then $J_n \subseteq J$, $J_n$ is a reduction of $I$ and $\ell(I/J_n) < \infty$.
Since $\cap_{n\geq c}J_n = L$, it follows that atleast one $J_m \neq J$. Put $K = J_m$.
\end{proof}

\section{Hilbert-Samuel Coefficients of $I^kM/I^{k+1}M$ as a function of $k$}

In this section $(R, \mm)$ is a Noetherian local ring, $I$ is an ideal in $R$, $M$  is a finitely generated $R$-module and $J$ is a reduction of $I$ with respect to $M$, i.e., $I^{m+1}M = JI^mM$ for some $m$.
If $\dim M = r$, then for $i = 0,1,\ldots, r$ let  $e_i(M)$ be the $i$th Hilbert coefficient of $M$ (with respect to $\mm$). In Corollary  \ref{locT}, we prove that $e_i(I^kM/I^{k+1}M)$,  $e_i(I^kM)$ and $e_i(I^k M/J^{k}M)$)
are polynomial functions  in $k$ for $k\gg 0$.

We might also ask similar questions for Hilbert coefficients with respect to an $\mm$-primary ideal in $R$. Our proofs for Question (1) and (3)  in the introduction also covers this slightly more  case. However for notational convenience and as a  matter of taste we stick to the maximal ideal.

As in Section 3, it is convenient to consider standard graded  $R$- algebras, $A = \bigoplus_{n \geq 0} A_n$ with $A_0 = R$,  and finitely generated graded $A$-modules. To ensure the existence of limit dimensions and a result regarding annihilators we show

\begin{Proposition}\label{annhilator}
 Let   $A = \bigoplus_{n \geq 0} A_n$ be a standard graded $R$-algebra with $A_0 = R$. Let $M = \bigoplus_{k \geq 0}M_k$ be a finitely generated $A$-module. Then
$\ann_R M_k =  \ann_R M_{k-1}$ for all $k \gg 0$. In particular $\dim M_k$ is constant for all $k \gg 0$.
\end{Proposition}

\begin{proof}
 Notice that $M_k = A_1 M_{k-1}$ for all $k \gg 0$. It follows that $ \ann_R M_{k-1} \subseteq \ann_R M_{k}$ for all
$k \gg 0$.

To prove the other inclusion we first assume that the residue field $K = R/\mm$ is infinite. Then by \cite[2.1,3.1]{t} there exists $u \in A_1$
which is $M$-filter regular. So the $R$-linear maps $M_{k-1} \rightarrow M_k$ given by multiplication by $u$ are injective
for $k \gg 0$. In particular $ \ann_R M_{k} \subseteq \ann_R M_{k-1}$ for  $k \gg 0$.  Thus the result holds when $K$ is infinite.

If $K$ is finite then we do the standard trick, i.e., we consider $R^{\prime} = R[x]_{\mm R[x]}$. The residue field of
$R^{\prime}$ is $K(x)$, which is infinite. Set $A^\prime = A \otimes_R R^{\prime}$
and $M_{k}^{\prime} = M \otimes_R R^{\prime}$ for $k \geq 0$. Notice that $M^{\prime} = M \otimes_R R^{\prime} = \bigoplus_{k \geq 0} M_{k}^{\prime}$ is a finitely generated $A^\prime$-module and $A^\prime$ is a finitely generated $R^\prime$-algebra. Thus as proved before one has
$\ann_{R^\prime} M_{k}^{\prime} =  \ann_{R^\prime} M_{k-1}^{\prime}$ for all $k \gg 0$.

As each $M_k$ is a finitely generated $R$-module we get
\[
 \ann_{R^\prime} M_{k}^{\prime}  = \ann_R M_{k} \otimes_R R^{\prime} \quad \text{for each} \  k \geq 1.
\]
In the first line of this proof we had shown that $ \ann_R M_{k-1} \subseteq \ann_R M_{k}$ for all
$k \gg 0$ without any assumptions on $K$. Therefore
\[
 \frac{\ann_R M_{k}}{ \ann_R M_{k-1}}\otimes_{R} R^{\prime}  =  \frac{\ann_{R^{\prime} } M_{k}^{\prime}}{ \ann_{R^{\prime} } M_{k-1}^{\prime}}  = 0 \quad \text{for all} \ k \gg 0.
\]
 $R^{\prime}$ is a faithfully flat $R$-algebra. So $\ann_R M_k =  \ann_R M_{k-1}$ for all $k \gg 0$.
\end{proof}

We call the $K$-algebra $A/\mm A$ the \emph{fiber} of $A$ and $\dim A/\mm A = s(A)$, the \emph{spread}
of $A$. These definitions are influenced by the case when  $A = \Rees(I)$. Then $A/\mm A = F(I)$ is called  the fiber-cone of $I$ and  $s(I) = \dim F(I) $ is  the analytic spread of $I$.

If $K = R/\mm$ is infinite,  then by the same argument as those that are used in the proof  of the
 existence of minimal reductions of an ideal, we can show that there exists $u_1,\ldots, u_{c} \in A_1$ (here $c= s(A)$) such that $A_k = (u_1,\ldots, u_{c})A_{k-1}$ for all
$k \gg 0$.

\begin{Remark}\label{base1}
\emph{Set $B = R[y_1,\ldots, y_c]$. Consider the ring homomorphism $\phi \colon B \rightarrow A$ which sends $y_i$ to
$u_i$ for each $i$. Then $A$, considered as a $B$-module via $\phi$, is a finitely generated $B$-module.}
 \end{Remark}

\begin{Theorem} \label{mtLocal}
With hypothesis as in Proposition~\ref{annhilator} let $d$ be the limit dimension of the family $\{ M_k \}_{k \geq 0}$. Then
 $e_i(M_k)$   is a polynomial of degree  $\leq d + s(A) - 1 -i$ for $k\gg 0$.
\end{Theorem}
\begin{proof}
We will prove the theorem in three steps.

\medskip
\noindent
\emph{Step I:} We first  consider the function $f(k,n) = \ell (M_k/\mm^{n+1}M_k)$ and show
\[
 \sum_{k,n \geq 0}f(k,n)z^nw^k = \frac{h(z,w)}{(1-z)^{p+1}(1-w)^{q}} \quad \text{where} \ h(z,w) \in \mathbb{Z}[z,w].
\]

 Consider the Rees ring $S = \bigoplus_{n \geq 0} (\mm A)^n$ and the associated graded ring
$T = \gr_{\mm A} A$ of the ideal $\mm A$ in $A$. Observe  that $(\mm A)^n = \bigoplus_{k \geq 0}\mm^nA_k$. So
$S = \bigoplus_{k,n \geq 0}\mm^nA_k$ is a naturally standard bigraded $R$-algebra with $S_{k,n} = \mm^nA_k$ and $S_{0,0} = R$.
It can also be easily checked that $T = \bigoplus_{k,n \geq 0}T_{k,n}$, where $T_{k,n} = \mm^nA_k/\mm^{n+1}A_k$
is a standard bigraded $K$-algebra.

Next consider the Rees module  $ E = \Rees(\mm A,M) = \bigoplus_{n \geq 0}(\mm A)^nM$ and the associated graded module
$N =\gr_{\mm A} M$ of $M$ with respect to the ideal $\mm A$. Clearly $E$ and  $N$ are  finitely generated bigraded $S$ and  $T$ modules,  respectively. Notice that $N = \bigoplus_{k,n \geq 0} \mm^n M_k/\mm^{n+1}M_k$.

By a bigraded version of the Hilbert-Serre theorem, cf. \cite[8.20]{MS} we get that
\[
 \sum_{k,n \geq 0} \ell\left(\frac{\mm^n M_k}{\mm^{n+1}M_k} \right)z^nw^k = \frac{h(z,w)}{(1-z)^{p}(1-w)^{q}} \quad \text{where} \ h(z,w) \in \mathbb{Z}[z,w].
\]
Therefore
\[
 \sum_{k,n \geq 0} \ell\left(\frac{ M_k}{\mm^{n+1}M_k} \right)z^nw^k = \frac{h(z,w)}{(1-z)^{p+1}(1-w)^{q}}.
\]

\medskip
\noindent
\emph{Step II:}  $e_i(M_k)$ is of  polynomial type in $k$  of degree $ \leq  p +q -1 -i$.

By Step 1. we get that  there exists   $k_0$ and $n_0 \gg 0 $,  such that
\[
 f(n,k) = \sum_{0 \leq i + j \leq p + q -1} a_{i,j}\binom{n+i}{i}\binom{k+j}{j} \quad \ \text{for all} \ n \geq n_0 \ \text{ and } \ k \geq k_0.
\]
$\text{Here} \ a_{i,j} \in \mathbb{Z}.$
Since $\dim M_k$ is $d$ for $k \gg 0$ we get that $a_{i,j} = 0$ for all $i > d$. Fix $k \geq k_0$. It follows that
 for all $n \gg 0$ we have
\[
 \ell\left(\frac{M_k}{\mm^{n+1}M_{k}} \right) = f(n,k) = \sum_{i = 0}^{d}\left(\sum_{0 \leq j \leq p +q -1-i}a_{i,j}\binom{k+j}{j} \right)\binom{n+i}{i}
\]
Hence
\[
 e_i(M_k) = \sum_{0 \leq j \leq p +q -1-i}a_{i,j}\binom{k+j}{j} \quad \text{for} \ i = 0, \ldots, d.
\]
This implies that $e_i(M_k)$ is a polynomial in $k$ of degree $ \leq p +q -1 -i$ for $k\gg 0$.

\medskip
\noindent
\emph{Step III:}  We show that $\displaystyle{ p \leq d  \quad \text{and} \quad q \leq s(A) }$.

 We first do several reduction steps.

\medskip
\noindent
\emph{Reduction step  (a):} We may assume $d = \dim M_k = \dim R$ for all $k \geq 0$.
In fact, by Proposition \ref{annhilator}, there exists $r$ such that $\ann_R M_k = \ann_R M_r = \qq$ (say) for $k \geq r$.
We consider the submodule $D = \bigoplus_{k \geq r} M_k$ of $M$. Notice that $D$ is a $A/\qq A$-module and
$s(A/\qq A) \leq s(A)$.  Thus (a) follows.

\medskip
\noindent
\emph{Reduction step  (b):} We may assume $K$ is infinite. Indeed, if this is not the case we may apply a suitable base field extension.

\medskip
\noindent
\emph{Reduction step  (c):} We may assume $A = R[y_1,\ldots, y_c]$, where  $c=s(A)$. This follows from Remark \ref{base1}.

\medskip
\noindent
Now let $J = (x_1,\ldots, x_d)$ be a minimal reduction of $\mm$. Then $JA$ is a reduction of $\mm A$. It follows that
$S = \bigoplus_{n \geq 0} (\mm A)^n$ is finitely generated as a $S^{\prime} = \bigoplus_{n \geq 0} (J A)^n$-module.
Therefore $N = E \otimes_R R/\mm$ is a finitely generated  $S^{\prime}\otimes_R R/\mm$-module.
Since $T^\prime = S^\prime \cong \Rees(J)[y_1,\ldots, y_c]$, we get $S^{\prime}\otimes_R R/\mm \cong  F(J)[y_1,\ldots, y_c]$.
(Here $\Rees(J)$ is the Rees algebra of $J$ and $F(J)$ is the fiber-cone of $J$). Since $J$ is a minimal reduction of $\mm$, it follows that $F(J)$ is a polynomial ring over $K$ in $d$-variables, say $t_1,\ldots,t_d$. Notice that
\begin{align}
 T^\prime &\cong k[t_1,\ldots,t_d, y_1,\ldots, y_c] \quad \text{where} \ \deg t_j = (0,1)  \ \text{and}  \ \deg y_i = (1,0) \label{base2-1}, \text{ and that}\\
 \ N &=  \bigoplus_{k,n \geq 0} \mm^n M_k/\mm^{n+1}M_k  \ \text{is a finitely generated, bigraded $T^\prime$-module.} \label{base2-2}
\end{align}
It follows that $p \leq d$ and $q \leq c \leq s(A)$, see \cite[8.20]{MS}. This completes Step III.

\medskip
Our result follows from Steps I,II and III.
\end{proof}

\begin{Corollary}\label{locT}
 Let $M$ be a finitely generated $R$-module. Let $I$ be an ideal in $R$ and let $J$ be a reduction of $I$ with
respect to $M$. Then for $k \gg 0$,
\begin{enumerate}[\rm (a)]
 \item
$e_i(I^kM/I^{k+1}M)$ is a polynomial in $k$ of degree $\leq \alpha_I(M) + s(I) - i -1$, where $0 \leq i \leq \alpha_I(M)$.
\item
$e_i(I^kM)$ is a polynomial in $k$ of degree $\leq \beta_I(M) + s(I) - i -1$, where $0 \leq i \leq \beta_I(M)$.
\item
  $e_i(I^k M/J^{k}M)$ polynomial in $k$ of degree $\leq \gamma^{I}_{J}(M) + s(J) - i -1$, where  $i = 0,\ldots, \gamma^{I}_{J}(M)$.
\end{enumerate}
\end{Corollary}
\begin{proof}
 The assertions (a) and (b) follow from Theorem \ref{mtLocal},
 since $ \bigoplus_{k \geq 0} I^kM/I^{k+1}M $ and $\bigoplus_{k \geq 0} I^kM  $ are finitely generated modules over $ \Rees(I)$.  The assertion (c) follows from Theorem \ref{mtLocal},
 since $ \bigoplus_{k \geq 0} I^kM / J^{k}M$ is a
 finitely generated modules over $ \Rees(J)$.
\end{proof}

The following result regarding multiplicities of modules over not-necessarily standard graded algebras is needed in Section 3.
\begin{Proposition}\label{XX}
Let $B = \bigoplus_{k \geq 0} B_k$ be a finitely generated $R$-algebra with $B_0 = R$. Let $M = \bigoplus_{k \geq 0}M_k$
be a finitely generated $B$-module. Then there exists $g \geq 1$ such that
\begin{enumerate}[\rm (a)]
 \item
$\dim M_{mg + i}$ is constant (say $d(i)$) for $m \gg 0$  for each $i = 0, 1,\ldots, g -1$.
\item
The function $f(k) = e_0(M_k)$ is  quasi-polynomial of period $g$ for $k\gg 0$.
\item
If $P_0, \ldots, P_{g-1}$ are  polynomials  such that for each $i$ we have  $P_i(m) = e_0(M_{mg + i})$ for $m \gg 0$,  then
$\deg P_i \leq \dim B/\mm B + d(i) -1$.
\end{enumerate}
\end{Proposition}
\begin{proof}
 We choose $g$ such that the Veronese subring $A = B^{(g)} $ is standard graded. Then $M ^{(g;i)} = \bigoplus_{k \geq 0}M_{kg +i}$ is a finitely generated $A$-module for each $i = 0,\ldots, g-1$.  Also $M = \bigoplus_{i =0}^{g-1}M ^{(g;i)}$. The assertions (a) and (b) follow from Theorem \ref{mtLocal}. To get the degree estimate notice that
\[
 \left(\frac{B}{\mm B}\right)^{(g)} = \frac{B^{(g)}}{\mm B^{(g)}} = \frac{A}{\mm A}.
\]
So $s(A) = \dim A/\mm A = \dim B/ \mm B$.
\end{proof}

\end{document}